\documentclass{article}
\usepackage{amsfonts}
\usepackage{amsmath}

\newtheorem{theorem}{Theorem}
\newtheorem{acknowledgement}[theorem]{Acknowledgement}

\newtheorem{corollary}[theorem]{Corollary}

\newtheorem{lemma}[theorem]{Lemma}

\newtheorem{proposition}[theorem]{Proposition}
\newtheorem{remark}[theorem]{Remark}

\input{tcilatex}
\begin{document}

\title{On cubic Berwald spaces}
\author{Nicoleta Brinzei \\
Transilvania University, Brasov, Romania}
\maketitle

\begin{abstract}
We show that, for Finsler spaces with cubic metric, Landsberg spaces are
Berwaldian. Also, for decomposable metrics, we determine specific conditions
for a space with cubic metric to be of Berwald type, thus refining the
result in \cite{Mat}.
\end{abstract}

\section{Introduction}

Spaces with cubic metric are studied by Matsumoto and Numata, \cite{Mat}, 
\cite{Mat2}. They are Finsler spaces in a wider sense, \cite{Rashewski}.

An interesting problem related to m-th root metric spaces is the following:
is any Landsberg space with m-th root metric Berwaldian?

A partial answer for spaces with cubic metric with fundamental function $F=%
\sqrt[3]{c_{1}\alpha ^{2}\beta +c_{2}\beta ^{3}}$ (where $\alpha ^{2}$ is a
pseudo-Riemannian metric and $\beta $ is a 1-form) is given by Lee and Jun, 
\cite{Lee}. In what follows, we generalize this result: namely, for all
cubic Finsler spaces $(M,F),$ $F=\sqrt[3]{a_{ijk}(x)y^{i}y^{j}y^{k}}$ with $%
a_{ijk}$ differentiable, if $(M,F)$ is of Landsberg type, then it is of
Berwald type.

Also, for spaces whose fundamental function is decomposable as a product of
two factors $\bar{F}^{3}=a\cdot b,$ between a Riemannian metric $a$ and a
1-form $b$ on $M,$ we show that $(M,\bar{F})$ is of Berwald type if and only
if the 1-form $b$ is parallelly transported with respect to the Levi-Civita
connection of $a.$ An analogous result is proven by Z. Shen for spaces with $%
(\alpha ,\beta )$-metrics of the form $F=\alpha \phi (\dfrac{\beta }{\alpha }%
),$ \cite{Shen1}.

The techniques we used mainly rely on expressing the involved geometrical
objects in terms of the third power $T=F^{3}$ of the fundamental function,
which is a polynomial function of the directional variables $y^{i}.$

\bigskip

\section{Spaces with cubic metric}

Let $M^{n}$ be a differentiable manifold of dimension $n$ and class $%
\mathcal{C}^{\infty },$ $TM$ its tangent bundle and $(x^{i},y^{i})$ the
coordinates in a local chart on $TM$. Let $F$ be the following function on $%
M,$ :%
\begin{equation}
F=\sqrt[3]{a_{ijk}(x)y^{i}y^{j}y^{k}}.  \label{1.1}
\end{equation}%
(with $a_{ijk}$ symmetric in all its indices) and%
\begin{equation}
T=F^{3}=a_{ijk}(x)y^{i}y^{j}y^{k}.  \label{1}
\end{equation}

In the following, for a function $f=f(x,y),$ we shall denote by " $,"$ and "~%
$\cdot $~" the partial derivatives w.r.t. $x$ and $y,$ respectively. Also,
if $N$ is a nonlinear connection on $TM,$ we denote by " ; " its associate
covariant derivative%
\begin{equation*}
f_{;l}=\frac{\delta f}{\delta x^{l}}=\frac{\partial f}{\partial x^{l}}%
-N_{~l}^{r}\frac{\partial f}{\partial y^{r}},\;f\in \mathcal{F}(TM)
\end{equation*}

and we denote by null index transvection by $y$ (for instance, $%
T_{i0}=T_{ij}y^{j})$.

\begin{remark}
\cite{projective} If $F=T^{1/m}$ is a Finslerian fundamental function on $M$%
, then the Hessian $[T_{ij}]$ is an invertible matrix, its inverse has the
entries:%
\begin{equation*}
T^{ij}=\dfrac{1}{m(m-1)F^{m-2}}\{(m-1)g^{ij}-(m-2)l^{i}l^{j}\},
\end{equation*}%
where $g^{ij}$ denotes the contravariant version of the usual Finslerian
metric tensor attached to $F$ and $l^{i}=\dfrac{y^{i}}{F}.$
\end{remark}

Hence, $T^{ij}$ and $T_{ij}$ can be used for raising and lowering indices of
tensors. Moreover, $T_{ij}$ are polynomial functions of $y,$ and $T^{ij}$
are rational functions of $y.$

\section{Geodesics and canonical spray}

In the following, we shall express the equations of geodesics of a cubic
metric space and the related geometric objects in terms of $T=F^{3}$ of the
fundamental function and of its derivatives.

\bigskip

Unit speed geodesics of $(M,F)$ are described by the Euler-Lagrange equation:%
\begin{equation*}
\dfrac{\partial F}{\partial x^{i}}-\dfrac{d}{dt}\left( \dfrac{\partial F}{%
\partial y^{i}}\right) =0.
\end{equation*}%
Taking into account the fact that, along such curves, $F(x,\dot{x})=1,$ the
above is equivalent to:%
\begin{equation*}
\dfrac{\partial T}{\partial x^{i}}-\dfrac{d}{dt}\left( \dfrac{\partial T}{%
\partial y^{i}}\right) =0.
\end{equation*}

An easy computation leads to:%
\begin{equation}
\dfrac{dy^{i}}{dt}+T^{ih}(T_{h,k}y^{k}-T_{,h})=0,~\ \ \ y^{i}=\dot{x}^{i}.
\label{geodesics}
\end{equation}

Consequently,

\begin{proposition}
\begin{enumerate}
\item In spaces with cubic metric the coefficients of the canonical spray, 
\cite{Anto}, \cite{Miron}, are rational functions of $(y^{i}),$ given by%
\begin{equation}
2G^{i}=T^{ih}(T_{h,k}y^{k}-T_{,h}).  \label{G}
\end{equation}

\item The canonical nonlinear connection has the coefficients:$%
N_{~j}^{i}=G_{~\cdot j}^{i}=\dfrac{1}{2}\{T_{~~\cdot
j}^{ih}(T_{h,k}y^{k}-T_{,h})+T^{ih}(T_{hj,k}y^{k}+T_{h,j}-T_{j,h})\}.$
\end{enumerate}
\end{proposition}

\bigskip

We denote in the following by $B\Gamma $ the Berwald connection, \cite{Anto}%
, \cite{Shen} determined by $F=\sqrt[3]{T}$ and by $G_{~jk}^{i}=G_{~\cdot
jk}^{i}$ its coefficients. According to (\ref{G}), for m-th root metric
spaces, $G_{~jk}^{i}$ are rational functions of $y.$

Also, let 
\begin{eqnarray*}
L_{~jk}^{i} &=&\dfrac{T^{ih}}{2}(\dfrac{\delta T_{hj}}{\delta x^{k}}+\dfrac{%
\delta T_{hk}}{\delta x^{j}}-\dfrac{\delta T_{jk}}{\delta x^{h}}), \\
T_{~jk}^{i} &=&\dfrac{T^{ih}}{2}(\dfrac{\partial T_{hj}}{\partial y^{k}}+%
\dfrac{\partial T_{hk}}{\partial y^{j}}-\dfrac{\partial T_{jk}}{\partial
y^{h}})=\dfrac{T^{ih}}{2}T_{hjk}.
\end{eqnarray*}%
denote the coefficients of the canonical metrical connection $C\Gamma $
attached to the Lagrange-type metric $T_{ij},$ \cite{Miron}.

\bigskip

\section{Specific Landsberg\&Berwald conditions for mth- root metrics}

\bigskip

There are a lot of alternative definitions of Landsberg and Berwald-type
Finsler spaces, \cite{Anto}, \cite{Dodson}. In the present paper, we shall
use the following:

A Finsler space $(M,F)$ is a \textit{Landsberg space} if: (1) the Cartan
tensor $C_{ijk}$ satisfies $C_{ijk|0}=0,$ where the covariant derivative is
taken with respect to the Berwald connection $B\Gamma ,$ or (2): the Berwald
connection $B\Gamma $ is metrical.

In Ladsberg spaces, the horizontal coefficients of the Cartan connection $%
F_{~jk}^{i}$ coincide with those of the Berwald connection: $%
F_{~jk}^{i}=G_{~jk}^{i}.$

A Finsler space is called a \textit{Berwald space} if: (1) with respect to $%
B\Gamma (N),$ there holds $C_{ijk|l}=0$ or (2) the coefficients $G_{~jk}^{i}$
of the Berwald connection are functions of $x^{i}$ alone: $%
G_{~jk}^{i}=G_{~jk}^{i}(x).$

The last statement is equivalent to the fact that the coefficients $G^{i}$
of the canonical spray are homogeneous polynomial functions of degree 2 in $%
y^{i}.$ There hold the inclusions:%
\begin{equation*}
Riemann~spaces\subset Berwald~spaces~\subset Landsberg~spaces.
\end{equation*}

For Finsler spaces with m-th root metric $(M,F)$, we get more convenient
such characterizations by using the third order derivatives $T_{ijk}$ (where 
$T=F^{m}$) instead of the Cartan tensor $C_{ijk}.$

Using the results in \cite{Shimada}, we have proven in \cite{projective},
that

\begin{proposition}
The horizontal coefficients $L_{~jk}^{i}$ of the canonical metrical
connection $C\Gamma \ $attached to the Hessian $T_{ij}$ coincide with those
of the Cartan connection of $(M,F).$ Hence, in Landsberg m-th root metric
spaces, we have $L_{~jk}^{i}=F_{~jk}^{i}=G_{~jk}^{i}.$
\end{proposition}

\begin{corollary}
An m-th root metric space $(M,F)$ is a Berwald space (resp. Landsberg space)
if and only if, w.r.t. the canonical metrical connection $C\Gamma (N),$ we
have $T_{ijk|l}=0$ (resp. $T_{ijk|0}=0$).
\end{corollary}

\bigskip

\section{Landsberg-Berwald equivalence}

In the following, we show that Landsberg spaces with cubic metrics are
Berwaldian.

Let 
\begin{equation*}
T=F^{3}=a_{ijk}(x)y^{i}y^{j}y^{k},
\end{equation*}%
with $a_{ijk}=a_{ijk}(x)$ of class at least 1, define a Landsberg space;
according to the results in the previous section, this means 
\begin{equation*}
T_{ijk|0}=0.
\end{equation*}

\bigskip

For a cubic metric, the third derivatives $T_{ijk}$ depend only on $x,$
which entails $\dfrac{\delta T_{ijk}}{\delta x^{l}}=\dfrac{\partial T_{ijk}}{%
\partial x^{l}}.$

Then,%
\begin{equation}
T_{ijk|l}=T_{ijk,l}-L_{~il}^{h}T_{hjk}-L_{~jl}^{h}T_{ihk}-L_{~kl}^{h}T_{ijh}.
\label{h-deriv1}
\end{equation}%
Taking into account that our space is a Landsberg one (i.e., $%
L_{~il}^{h}=G_{~il}^{h}$ etc.), we have%
\begin{equation*}
T_{ijk|0}=T_{ijk,l}y^{l}-N_{~i}^{h}T_{hjk}-N_{~j}^{h}T_{ihk}-N_{~k}^{h}T_{ijh}=0.
\end{equation*}%
Deriving by $y^{l}$ and taking into account that $T_{ijk}$ depend only on $%
x, $ we get%
\begin{equation*}
T_{ijk,l}-L_{~il}^{h}T_{hjk}-L_{~jl}^{h}T_{ihk}-L_{~kl}^{h}T_{ijh}=0,
\end{equation*}%
which is nothing but $T_{ijk|l}=0.$ We have thus obtained

\begin{proposition}
Let $(M,F)$ be a space with cubic metric $F=\sqrt[3]{%
a_{ijk}(x)y^{i}y^{j}y^{k}}.$ If the functions $a_{ijk}$ are of class at
least one, then there holds the implication:
\end{proposition}

\begin{center}
$(M,F)$ is a Landsberg space $\Rightarrow $ $(M,F)$ is a Berwald space.
\end{center}

\bigskip

Further, for spaces with cubic metric, the inclusion Riemannian spaces $%
\subset $ Berwald spaces is strict. Namely, the Berwald-Moor conformal space
with%
\begin{equation*}
T=F^{3}=e^{\sigma (x)}y^{1}y^{2}y^{3},
\end{equation*}%
where $\sigma (x)$ is a differentiable function, provides an example of
Berwald cubic space, which is non-Riemannian.

\bigskip

\section{Decomposable cubic metrics}

Let us consider a space $(M,F=\sqrt[3]{T}),$ where $T$ decomposes as a
product%
\begin{equation}
T=a\cdot b  \label{decomp}
\end{equation}%
where $a=\gamma _{ij}(x)y^{i}y^{j}$ is a Riemannian metric and $b=b_{i}(x)$
is a 1-form, such that:

\begin{equation*}
\left\Vert b\right\Vert ^{2}=\gamma ^{ij}b_{i}b_{j}=1.
\end{equation*}

For cubic spaces with $T=F^{3}$ as in (\ref{decomp}), it is proven in \cite%
{Mat} that the space is a Berwald one if and only if there exists some
1-form $f\in \mathcal{X}^{\ast }M$ such that%
\begin{equation*}
\gamma _{ij|k}=f_{k}(x)\gamma _{ij};~~\ \ b_{i|k}=-f_{k}(x)b_{i},
\end{equation*}%
where the covariant derivative is taken with respect to the Berwald
connection determined by the "whole" fundamental function $F=\sqrt[3]{ab}.$

In the following, we shall find the relation between $a$ and $b$ such that
the space $(M,F=\sqrt[3]{ab})$ is Berwaldian; more precisely, we shall take
into consideration the covariant derivatives%
\begin{equation*}
\nabla _{i}b_{j},
\end{equation*}%
where $\nabla $ denotes the Levi-Civita connection attached to $\gamma
_{ij}. $

By direct computation, we get

\begin{lemma}
If $a=\gamma _{ij}(x)y^{i}y^{j}$ is a Riemannian metric and $b=b_{i}(x)$ is
a 1-form with $\gamma ^{ij}b_{i}b_{j}=1,$ then:

\begin{enumerate}
\item The Hessian matrix $[T_{ij}]$ is invertible iff%
\begin{equation*}
\Delta :=4b^{2}-a
\end{equation*}%
does not vanish.

\item The inverse matrix has the entries%
\begin{equation}
T^{ij}=\dfrac{1}{2b\Delta }(\Delta \gamma
^{ij}-2bb^{i}y^{j}-2bb^{j}y^{i}+ab^{i}b^{j}+y^{i}y^{j}),  \label{inverse}
\end{equation}%
where the indices of $b$ were raised by $\gamma ^{ih}:$ $b^{i}=\gamma
^{ih}b_{h}.$
\end{enumerate}
\end{lemma}

Furhter, in \cite{Anto}, p. 110-111, it is proven the following result:

\begin{lemma}
, \cite{Anto}: If $(M,F)$ and $(M,\bar{F})$ are two Finsler spaces on the
same underlying manifold, then the local coefficients of the corresponding
canonical sprays are related by%
\begin{equation}
2\bar{G}^{i}=2G^{i}+\dfrac{\bar{F}_{|0}y^{i}}{\bar{F}}-\bar{F}\bar{g}%
^{ij}r_{j}(\bar{F}),  \label{difference}
\end{equation}%
where $|$ denotes Berwald covariant derivative determined by $F$ and%
\begin{equation*}
r_{j}(S)=S_{|j}-y^{r}S_{|r\cdot j},~\ \ \forall S\in \mathcal{F}(TM).
\end{equation*}
\end{lemma}

In the following, we shall express the above in terms of the m-th power of $%
\bar{F},$ $m\geq 2;$ hence, let for the moment 
\begin{equation*}
T=\bar{F}^{m}.
\end{equation*}
Then, there hold the relations:

\begin{itemize}
\item 
\begin{equation}
\dfrac{\bar{F}_{|0}y^{i}}{\bar{F}}=\dfrac{1}{m}\dfrac{T_{|0}y^{i}}{T}.
\label{*}
\end{equation}

\item The contravariant Finslerian metric tensor $\bar{g}^{ij}$ is expressed
in terms of $T$ as%
\begin{equation*}
\bar{g}^{ij}=\dfrac{T^{~-\tfrac{2}{m}}}{m-1}\left(
Tm(m-1)T^{ij}+(m-2)y^{i}y^{j}\right) .
\end{equation*}

\item $r_{j}(\bar{F})=\dfrac{1}{m}T^{\tfrac{1}{m}-2}\left( Tr_{j}(T)+\dfrac{%
m-1}{m}T_{j}T_{|0}\right) ;$

\item $y^{j}r_{j}(T)=(1-m)T_{|0}.$
\end{itemize}

Then, the last term in (\ref{difference}) is%
\begin{equation*}
\begin{array}{l}
\bar{F}\bar{g}^{ij}r_{j}(\bar{F})=T^{\tfrac{1}{m}}\dfrac{T^{~-\tfrac{2}{m}}}{%
m-1}\left( Tm(m-1)T^{ij}+(m-2)y^{i}y^{j}\right) \cdot \\ 
\cdot \dfrac{1}{m}T^{\tfrac{1}{m}-2}\left( Tr_{j}(T)+\dfrac{m-1}{m}%
T_{j}T_{|0}\right) = \\ 
=\dfrac{T^{-2}}{m(m-1)}%
\{T^{2}m(m-1)T^{ij}r_{j}(T)+(m-2)y^{i}y^{j}Tr_{j}(T)+(m-1)^{2}TT^{ij}T_{j}T_{|0}+
\\ 
+\dfrac{(m-2)(m-1)}{m}y^{i}y^{j}T_{j}T_{|0}\}=T^{ij}r_{j}(T)-\dfrac{m-2}{m}%
T^{-1}y^{i}T_{|0}+\dfrac{1}{m}T^{-1}T_{|0}y^{i}+ \\ 
+\dfrac{m-2}{m}T^{-1}y^{i}T_{|0}=T^{ij}r_{j}(T)+\dfrac{1}{m}%
T^{-1}T_{|0}y^{i}.%
\end{array}%
\end{equation*}%
Replacing into (\ref{difference}) and taking (\ref{*}) into account, we get

\begin{lemma}
If $(M,F)$ and $(M,\bar{F})$ are two Finsler spaces on the same underlying
manifold, then the coefficients of the corresponding canonical sprays are
related by%
\begin{equation}
2\bar{G}^{i}=2G^{i}-T^{ij}r_{j}(T),
\end{equation}%
where $|$ denotes Berwald covariant derivative determined by $F$ and%
\begin{equation*}
T=\bar{F}^{m},~m\geq 2,\ \ \ r_{j}(T)=T_{|j}-y^{r}T_{|r\cdot j}.
\end{equation*}
\end{lemma}

\bigskip

\bigskip

We shall also use the following relations, which can be deduced by direct
computation:%
\begin{eqnarray}
&&r_{j}(b)=(\nabla _{j}b_{r}-\nabla _{r}b_{j})y^{r};  \notag \\
&&y^{j}r_{j}(b)=0;  \label{aux} \\
&&T^{ij}b_{j}=\dfrac{1}{2\Delta }(2bb^{i}-y^{i});~T^{ij}a_{\cdot j}=\dfrac{1%
}{\Delta }(2by^{i}-b^{i}a);  \notag \\
&&\left\Vert b\right\Vert =1\Rightarrow b^{i}\nabla _{j}b_{i}=0.  \notag
\end{eqnarray}

Let now $G^{i}$ be determined by the Riemannian metric $\gamma _{ij}(x),$
where $a=\gamma _{ij}(x)y^{i}y^{j}$, and $\bar{G}^{i}$, by $T=\bar{F}%
^{3}=a\cdot b$ as above. Then, $|_{i}=\nabla _{i},$ and%
\begin{equation*}
r_{j}(T)=\nabla _{j}(ab)-y^{r}\dfrac{\partial }{\partial y^{j}}\nabla
_{r}(ab),
\end{equation*}%
and taking into account that 
\begin{equation*}
\nabla _{j}a=0,
\end{equation*}%
we get%
\begin{equation*}
r_{j}(T)=ar_{j}(b)-a_{\cdot j}\nabla _{0}b,
\end{equation*}%
where $\nabla _{0}b=y^{r}\nabla _{r}b.$

\bigskip

The cubic space $(M,\bar{F})$ is a Berwald one if and only if the functions $%
2\bar{G}^{i}$ are polynomial in $y^{i}.$ This is equivalent to the fact that
the difference%
\begin{equation*}
2B^{i}:=2\bar{G}^{i}-2G^{i}=-T^{ij}r_{j}(T)
\end{equation*}%
is a polynomial function of degree 2 in $y.$ There holds

\bigskip

\begin{theorem}
The space $(M,F=\sqrt[3]{T}),$ where $T$ decomposes as a product%
\begin{equation}
T=a\cdot b
\end{equation}%
where $a=\gamma _{ij}(x)y^{i}y^{j}$ is a Riemannian metric and $b=b_{i}(x)$
is a 1-form, such that:%
\begin{equation*}
\left\Vert b\right\Vert ^{2}=\gamma ^{ij}b_{i}b_{j}=1
\end{equation*}

\begin{enumerate}
\item is of Berwald type, if and only if $b$ is parallel with respect to $a:$%
\begin{equation*}
\nabla _{i}b_{j}=0,~\ \ \forall i,j=1,...,n.
\end{equation*}
\end{enumerate}
\end{theorem}

\textbf{Proof:}

Let us suppose that $(M,\bar{F}=\sqrt[3]{ab})$ is Berwaldian and let us fix
some arbitrary $x\in M$. Then $2B^{i}$ are polynomials of degree 2 and
hence, so are $2B^{i}b_{i}.$ By (\ref{aux}), we have $T^{ij}b_{j}=\dfrac{1}{%
2\Delta }(2bb^{i}-y^{i})$, consequently,%
\begin{eqnarray*}
-2B^{i}b_{i} &=&\dfrac{1}{2\Delta }(2bb^{j}-by^{j})r_{j}(T)=\dfrac{1}{%
2\Delta }(2bb^{j}-y^{j})(ar_{j}(b)-a_{\cdot j}\nabla _{0}b)= \\
&=&\dfrac{1}{\Delta }(abb^{j}r_{j}(b)-2b^{2}\nabla _{0}b+a\nabla _{0}b).
\end{eqnarray*}%
But, $a-2b^{2}=2b^{2}-\Delta ,$ so we can write%
\begin{equation*}
-2B^{i}b_{i}=\dfrac{1}{\Delta }\{abb^{j}r_{j}(b)+(2b^{2}-\Delta )\nabla
_{0}b\}=-\nabla _{0}b+\dfrac{1}{\Delta }\{abb^{j}r_{j}(b)+2b^{2}\nabla
_{0}b\}.
\end{equation*}%
Since the latter is a polynomial, $\Delta $ divides the polynomial $%
abb^{j}r_{j}(b)+2b^{2}\nabla _{0}b=b(ab^{j}r_{j}(b)+2b\nabla _{0}b).$ Since $%
a$ does not decompose in factors, $a$ and $b$ have no common factors; we
notice that, in this case, $b$ and $\Delta $ are also relatively prime, hence%
\begin{equation*}
\Delta ~|~ab^{j}r_{j}(b)+2b\nabla _{0}b.
\end{equation*}%
Again, we have $a=4b^{2}-\Delta ,$ and we get that $\Delta
~|~4b^{2}b^{j}r_{j}(b)+2b\nabla _{0}b=2b(2bb^{j}r_{j}(b)+\nabla _{0}b),$that
is, 
\begin{equation*}
\Delta ~|~(2bb^{j}r_{j}(b)+\nabla _{0}b).
\end{equation*}%
Both hand sides of the above are polynomials of degree 2 in $y^{i},$ hence
there exists some $f=f(x)$ such that%
\begin{equation}
(2bb^{j}r_{j}(b)+\nabla _{0}b)=f(x)\Delta .  \label{**}
\end{equation}

By identifiying the coefficients in the above relation and taking into
account that, by (\ref{aux}) $b^{i}\nabla _{j}b_{i}=0$, we get%
\begin{equation*}
2b_{i}b^{j}\nabla _{j}b_{r}+2b_{r}b^{j}\nabla _{j}b_{i}+\nabla
_{r}b_{i}+\nabla _{i}b_{r}=f(x)(8b_{i}b_{r}-2\gamma _{ir}).
\end{equation*}%
Contracting with $b^{i}$ and taking into account that $b^{i}b_{i}=1,$ the
above leads to%
\begin{equation}
b^{i}\nabla _{i}b_{r}=2b_{r}f(x),  \label{term1}
\end{equation}%
which yields%
\begin{equation}
b^{j}r_{j}(b)=b^{j}(\nabla _{j}b_{r}-\nabla _{r}b_{j})y^{r}=b^{j}\nabla
_{j}b_{0}=2bf(x).  \label{br}
\end{equation}

Replacing into (\ref{**}), we have $4b^{2}f(x)+\nabla _{0}b=f(x)\Delta
=f(x)(4b^{2}-a);$ we obtain that%
\begin{equation}
\nabla _{0}b=-af(x).  \label{term2}
\end{equation}

\bigskip

Let us come back now to the expression of $2B^{i}:$%
\begin{equation*}
-2B^{i}=T^{ij}(ar_{j}(b)-a_{\cdot j}\nabla _{0}b)
\end{equation*}

The last term, $T^{ij}a_{\cdot j}\nabla _{0}b$ is%
\begin{equation*}
T^{ij}a_{\cdot j}\nabla _{0}b=\dfrac{1}{\Delta }(2by^{i}-b^{i}a)\nabla _{0}b=%
\dfrac{-a}{\Delta }(2by^{i}-b^{i}a)f(x).
\end{equation*}%
The first one, $T^{ij}ar_{j}(b),$ is%
\begin{eqnarray*}
T^{ij}ar_{j}(b) &=&\dfrac{a}{2b\Delta }(\Delta \gamma
^{ij}-2bb^{i}y^{j}-2bb^{j}y^{i}+ab^{i}b^{j}+y^{i}y^{j})r_{j}(b)= \\
&=&\dfrac{a}{2b\Delta }(\Delta \gamma
^{ij}r_{j}(b)-0-4b^{2}y^{i}f(x)+2abb^{i}f(x)+0).
\end{eqnarray*}%
Then, 
\begin{equation*}
-2B^{i}=\dfrac{a}{2b\Delta }\{\Delta \gamma
^{ij}r_{j}(b)-4b^{2}y^{i}f(x)+2abb^{i}f(x)\}+\dfrac{2ab}{2b\Delta }%
(2by^{i}-b^{i}a)f(x)
\end{equation*}

The commom denominator $2b\Delta $ has to divide the numerator. In
particular, $b$ has to divide the numerator. The only term which does not
contain $b$ explicitely as a factor is 
\begin{equation*}
a\Delta \gamma ^{ij}r_{j}(b).
\end{equation*}%
Since $b$ has no common factors neither with $a$, nor with $\Delta ,$ $b$
has to divide the polynomial $\gamma ^{ij}r_{j}(b)$ (of degree 1). That is,
there exists some $\phi =\phi (x)$ such that $\gamma ^{ij}r_{j}(b)=\phi
^{i}(x)b.$ Lowering the indices, 
\begin{equation*}
r_{j}(b)=\phi _{j}(x)b.
\end{equation*}%
But, since $y^{j}r_{j}(b)=0,$ we get $0=y^{j}r_{j}(b)=(y^{j}\phi _{j}(x))b.$%
Together with $b\not=0,$ this yields $y^{j}\phi _{j}(x)=0,$ or%
\begin{equation*}
\phi _{j}=0,
\end{equation*}%
which is nothing but $r_{j}(b)=0.$ The latter means actually%
\begin{equation}
\nabla _{r}b_{i}-\nabla _{i}b_{r}=0.  \label{***}
\end{equation}

Let's now look at relation (\ref{term1}):%
\begin{equation}
b^{i}\nabla _{i}b_{r}=3b_{r}f(x),
\end{equation}

By (\ref{***}), it is equivalent to 
\begin{equation*}
b^{i}\nabla _{r}b_{i}=3b_{r}f(x).
\end{equation*}%
According to (\ref{aux}), we have $b^{i}\nabla _{r}b_{i}=0;$ the left hand
side of the above is 0, hence%
\begin{equation*}
f(x)=0,
\end{equation*}%
which yields, together with (\ref{term2}),%
\begin{equation*}
\nabla _{0}b=\nabla _{j}b_{i}y^{i}y^{j}=0.
\end{equation*}%
The latter, together with (\ref{***}), leads to%
\begin{equation*}
\nabla _{r}b_{i}=0,
\end{equation*}%
q.e.d.

The converse statement is obvious.

\begin{remark}
If $(M,\bar{F})$ is of Berwald type, then 
\begin{equation*}
2B^{i}:=2\bar{G}^{i}-2G^{i}=-T^{ij}r_{j}(T)=0,
\end{equation*}%
consequently, it has the same geodesics as the Riemannian space $(M,a=\gamma
_{ij}(x)y^{i}y^{j}).$
\end{remark}

\bigskip

\begin{acknowledgement}
Special thanks to prof. V. Balan for the helpful discussions and advices.

The work was supported by the Romanian Academy grant No.5 / 5.02.2008.
\end{acknowledgement}


\bigskip

\bigskip

\begin{thebibliography}{99}
\bibitem{Anto} Antonelli, P.L, Ingarden, R.S., Matsumoto, M., \textit{The
Theory of Sprays and Finsler Spaces with Applications in Physics and
Biology, }Kluwer Acad. Publ., 1993.

\bibitem{Shen} Bao, D., Chern, S.S., Shen, Z, \textit{An Introduction to
Riemann-Finsler Ge\-o\-me\-try} (Graduate Texts in Mathematics; 200),
Springer Verlag, 2000.

\bibitem{projective} Brinzei, N., \textit{Projective relations for m-th root
metric spaces}, arXiv:0711.4781v1.

\bibitem{Dodson} Dodson, C.T.J., \textit{A short review on Landsberg spaces}%
, http://www.maths.manchester.ac.uk/\symbol{126}%
kd/PREPRINTS/RevLandsberg.pdf.

\bibitem{Lee} Lee, Il-Yong; Jun, Dong-Gum, \textit{On two-dimensional
Landsberg space of a cubic Finsler space,} East Asian Math. J. 19, No.2,
305-316 (2003). [ISSN 1226-6973]

\bibitem{Mat} Matsumoto, M., Numata, S., \textit{On Finsler Spaces With a
Cubic Metric,} Tensor, N.S., 33(1979), 153-162.

\bibitem{Mat2} Matsumoto, M., \textit{Theory of Finsler spaces with m-th
root metric II, Publ. Math. Debrecen, 49(1996), 135-155.}

\bibitem{Miron} Miron, R., Anastasiei, M., \textit{Vector Bundles. Lagrange
Spaces. Applications in Relativity Theory}, (in Romanian), Ed. Acad.,
Bucharest, 1987.

\bibitem{Rashewski} Rashevsky, P. K. : \textit{The Geometrical theory of
partial differential equations}, Second Edition, Editorial USSR, M. 2003 (in
Russian).

\bibitem{Shimada} Shimada, H., \textit{On Finsler Spaces with Metric }$L=%
\sqrt[m]{a_{i_{1}i_{2}...i_{m}}y^{i_{1}}y^{i_{2}}...y^{i_{m}}},$ Tensor,
N.S., 33(1979), 365-372.

\bibitem{Shen1} Shen, Z. \textit{On Landsberg }$(\alpha ,\beta )$\textit{\
metrics}, http://www.math.iupui.edu/\symbol{126}zshen/Research/papers/

LandsbergCurvatureAlphaBeta2006.pdf.
\end{thebibliography}
\end{document}